\documentclass[11pt]{article}
\usepackage{amsfonts}
\usepackage{latexsym}

\newcommand{\bbz}{{\mathbf Z}}
\newcommand{\bbn}{{\mathbf N}}
\newcommand{\bbp}{{\mathbf P}}
\newcommand{\qed}{\hfill$\Box$}
\newcommand{\inv}{{\it inv}}
\newcommand{\Des}{{\it Des}}
\newcommand{\des}{{\it des}} 
\newcommand{\maj}{{\it maj}}
\newcommand{\Neg}{{\it Neg}}
\newcommand{\fmaj}{{\it fmaj}}

\newcommand{\sg}{{\sigma}}

\newtheorem{thm}{Theorem}[section]
\newtheorem{pro}[thm]{Proposition}
\newtheorem{lem}[thm]{Lemma}

\newtheorem{con}[thm]{Conjecture}
\newtheorem{cor}[thm]{Corollary}
\newtheorem{fac}[thm]{Fact}
\newtheorem{obs}[thm]{Observation}
\newtheorem{prb}[thm]{Problem}

\newtheorem{df}[thm]{Definition}
\newtheorem{rem}[thm]{Remark}

\newcommand{\eqdef}{:=}

\begin{document}
\pagestyle{myheadings}

\title{Equi-distribution over Descent Classes\\ of the Hyperoctahedral Group}
\bibliographystyle{acm}
\author{Ron M.\ Adin%
\thanks{Department of Mathematics and Statistics, Bar-Ilan University,
Ramat-Gan 52900, Israel. 
Email: {\tt radin@math.biu.ac.il} }\ $^\S$  
\and Francesco Brenti%
\thanks{Dipartimento di Matematica,
Universit\'{a} di Roma ``Tor Vergata'',
Via della Ricerca Scientifica,
00133 Roma, Italy. Email: {\tt brenti@mat.uniroma2.it}}\ $^\S$ 
\and Yuval Roichman%
\thanks{Department of Mathematics and Statistics, Bar-Ilan University,
Ramat-Gan 52900, Israel. 
Email: {\tt yuvalr@math.biu.ac.il} } 
\thanks{The research of all authors was partially supported by 
the EC's IHRP Programme, within the Research Training Network 
``Algebraic Combinatorics in Europe'', grant HPRN-CT-2001-00272}}
\date{submitted: June 9, 2004; revised: August 24, 2005}

\maketitle

\begin{abstract}
A classical result of MacMahon shows that the length function and 
the major index are equi-distributed over the symmetric group.
Foata and Sch\"utzenberger gave a remarkable refinement and proved that
these parameters are equi-distributed over inverse descent classes,
implying bivariate equi-distribution identities.
Type $B$ analogues of these results, refinements and consequences are given 
in this paper.
\end{abstract}

\section{Introduction}


Many combinatorial identities on groups 
are motivated by the fundamental works of MacMahon~\cite{MM}. 
Let $S_n$ be the symmetric group acting on $1,\dots,n$. 
We are interested in a refined enumeration of permutations according to 
(non-negative, integer valued) combinatorial parameters.
Two parameters that have the same generating function are said to be 
{\em equi-distributed}.
MacMahon~\cite{MM} has shown, about a hundred years ago, that  
the inversion number and the major index
statistics are equi-distributed on $S_n$ (Theorem~\ref{t.MM} below).
In the last three decades MacMahon's theorem has received
far-reaching refinements and generalizations.
Bivariate distributions were first studied by Carlitz~\cite{Ca}.
Foata~\cite{F} gave a bijective proof of MacMahon's theorem; 
then Foata and Sch\"utzenberger~\cite{FS} applied this bijection 
to refine MacMahon's identity, proving that
the inversion number and the major index are equi-distributed over 
subsets of $S_n$ with prescribed descent set of the inverse permutation 
(Theorem~\ref{FS} below).
Garsia and Gessel~\cite{GaG} extended the analysis to multivariate 
distributions. In particular, they gave an independent proof of the 
Foata-Sch\"utzenberger theorem, relying on an
explicit and simple generating function (see Theorem~\ref{GG} below).
Further refinements and analogues of the Foata-Sch\"utzenberger 
theorem were later found, involving left-to-right minima and 
maxima~\cite{BW} 
and pattern-avoiding permutations~\cite{RR, AR321}.
For a representation theoretic application of Theorem~\ref{FS}
see~\cite{R-Disc}.

\bigskip

Since the length and descent parameters may be defined via
the Coxeter structure of the symmetric group,
it is very natural to look for analogues of the above theorems
in other Coxeter groups.
This is a challenging open problem.
In this paper we focus on the hyperoctahedral group $B_n$, 
namely the classical Weyl group of type $B$.
Our goal is to find a type $B$ analogue of the
Foata-Sch\"utzenberger theorem (Theorem~\ref{FS}). To solve this we 
have to choose an appropriate type $B$ extension of the major index
among the many candidates, which were introduced and studied 
in~\cite{CF94, CF95a, CF95b, Rei93a, Rei93b, Ste94, FK}.
It turns out that the {\it flag-major index}, which was introduced in~\cite{AR} 
and further studied in~\cite{ABR, ABR2, HLR, AGR},
has the desired property:
the flag major index and the length function 
are equi-distributed on inverse descent classes of $B_n$.
In fact, we obtain a slight 
refinement of this result, involving the ``last digit'' parameter.
This parameter is involved in several closely related identities on $S_n$, 
see e.g.~\cite{AR321, AGR, RR}.
Our refinement also implies a MacMahon type theorem for the classical
Weyl group of type $D$, which has recently been proved in~\cite{BC}.
A summary of the results of this paper appeared in~\cite{ABR3}.

\bigskip

The rest of the paper is organized as follows.
Definitions, notation and necessary background are given in Section 2.
The main results are listed in Section 3.
Proofs of the main theorems are given in Section 4.
Section 5 contains open problems, an application to Weyl groups of type $D$
and remarks regarding different versions of the flag major index.
Finally, in the appendix (Section 6), we give an alternative proof of 
our type $B$ analogue of the Foata-Sch\"utzenberger theorem.

\section{Background and Notation}

\subsection{Notation}

Let $\bbp := \{ 1,2,3,\ldots \}$,
$\bbn := \bbp \cup \{0\} $, 
and $\bbz$ the ring of integers.
For $n \in \bbp$ let $[n] := \{ 1,2,\ldots,n \}$, and also $[0] := \emptyset$. 
Given $m, n\in \bbz$, $m \leq n$, let $[m,n] := \{ m, m+1,\ldots ,n \}$. 
For $n\in \bbp$ denote $[\pm n] := [-n,n] \setminus \{0\}$.
For $S \subset \bbn$ write $S = \{ a_1,\ldots,a_r \} _{<}$ to mean that 
$S= \{ a_1,\ldots,a_r \}$ and $a_1 < \ldots < a_r$.
The cardinality of a set $A$ will be denoted by $|A|$.

For $n\in \bbn$ denote 
\begin{eqnarray*}
[n]_q &:=& \frac{1-q^n}{1-q};\cr
[n]_q! &:=& \prod\limits_{i=1}^n \,[i]_q \qquad(n\ge1), \qquad [0]_q! := 1.
\end{eqnarray*}

For $n_1,\ldots,n_t\in \bbn$ such that $n_1 + \ldots + n_t =n$ define
the {\em $q$-multinomial coefficient}
$$ 
\left[n \atop n_1,\ldots,n_t \right]_q := \frac{[n]_q!}{[n_1]_q!\ldots[n_t]_q!}
$$
and use a shorter notation for the {\em $q$-binomial coefficient}
$$
{n\brack k}_q := {n\brack k,\,n-k}_q\qquad(0\le k\le n).
$$

\medskip

Given a statement $P$ we will sometimes find it convenient to let
\[
\chi (P) := \left\{ \begin{array}{ll}
1, &  \mbox{if $P$ is true,} \\
0, &  \mbox{if $P$ is false.}
\end{array} \right. 
\]

Given a sequence $\sigma = (a_{1}, \ldots ,a_{n}) \in \bbz^n$ 
we say that a pair $(i,j) \in [n] \times [n]$ is an 
{\em inversion} of $\sigma $ if $i<j$ and $a_{i}>a_{j}$.
We say that $i \in [n-1]$ is a {\em descent} of $\sigma $ if
$a_{i}>a_{i+1}$. We denote by $\inv(\sigma )$ (respectively, 
$\des(\sigma )$) the number of inversions (respectively, descents) of
$\sigma $. We also let
\[ 
maj (\sigma ) := \sum _{ \{ i : \; a_{i} > a_{i+1} \}}  i 
\]
and call it the {\em major index} of $\sigma$.

\smallskip

For $M=\{m_1,\dots,m_t\}_<\subseteq [0,n-1]$
denote $m_0:=0$, $m_{t+1}:=n$, and let
$$
{n \brack \Delta M}_q := {n \brack m_1-m_0,m_2-m_1,\dots,m_{t+1}-m_t}_q.
$$

Let $S_1,\ldots,S_k$ be sequences of distinct integers which are 
pairwise disjoint as sets. A sequence $S$ is a {\em shuffle} of $S_1,\ldots,S_k$
if $S$ is a disjoint union of $S_1,\ldots,S_k$ (as sets) and 
the elements of each $S_i$ appear in $S$ in the same order as in $S_i$.

\subsection{The Symmetric Group}

Let $S_n$ be the symmetric group on $[n]$.
Recall that $S_n$ is a Coxeter group with respect to the set of 
Coxeter generators $S:=\{s_i\,|\,1\le i\le n-1\}$, 
where $s_i$ may be interpreted as
the adjacent transposition $(i,i+1)$. 

If $\pi\in S_n$ then the classical combinatorial statistics 
(defined in the previous subsection)
of the sequence $(\pi(1),\dots,\pi(n))$ 
may also be defined via the Coxeter generators:
%
the inversion number $inv(\pi)$ is equal to the {\em length} $\ell(\pi)$
of $\pi$ with respect to the set of generators $S$; 
the {\em descent set} of $\pi$ is
$$
\Des(\pi):=
\{1\le i< n \,|\, \pi(i) > \pi(i+1)\} =
\{1\le i< n \,|\, \ell(\pi)>\ell(\pi s_i)\};
$$
the {\it descent number} of $\pi$ is
$\des(\pi):=|\Des(\pi)|$;
and the {\em major index} of $\pi$ is the sum (possibly zero)
$$ 
\maj(\pi) := \sum\limits_{i\in \Des(\pi)} i.
$$
The {\it inverse descent class} in $S_n$ corresponding to $M \subseteq [n-1]$
is the set $\{\pi\in S_n \,|\, \Des(\pi^{-1})= M\}$.
Note the following relation between inverse descent classes and shuffles. 

\begin{obs}\label{ObsA}
Let $\pi\in S_n$ and $M = \{m_1,\ldots,m_t\}_< \subseteq [n-1]$.
Then: $\Des(\pi^{-1})\subseteq M$ if and only if 
$(\pi(1),\ldots,\pi(n))$ is a shuffle of the following increasing sequences:
$$
(1,\ldots,m_1),
$$
$$
(m_1+1,\ldots,m_2),
$$
$$
\vdots
$$
$$
(m_t+1,\ldots,n).
$$
\end{obs}

\medskip

MacMahon's classical theorem asserts that the length function and 
the major index are equi-distributed on $S_n$.

\begin{thm}\label{t.MM}{\rm (MacMahon's Theorem~\cite{MM})}
$$
\sum_{\pi\in S_n} q^{\ell(\pi)}=
\sum_{\pi\in S_n} q^{\maj(\pi)} =[n]_q ! .
$$
\end{thm}

Foata~\cite{F} gave a bijective proof of this theorem.
Foata and Sch\"utzen\-berger~\cite{FS} applied this bijection to prove the following 
refinement.


\begin{thm}\label{FS}
{\rm (The Foata-Sch\"utzenberger Theorem~\cite[Theorem 1]{FS})}\\
For every subset $M\subseteq [n-1]$,
$$
\sum_{\{\pi\in S_n\,|\,\Des(\pi^{-1})=M\}} q^{\ell(\pi)}=
\sum_{\{\pi\in S_n\,|\,\Des(\pi^{-1})=M\}} q^{\maj(\pi)}.
$$
\end{thm}

\medskip

This theorem implies

\begin{cor}\label{FS1} 
$$
\sum\limits_{\pi\in S_n} q^{\ell(\pi)} t^{des(\pi^{-1})}=
\sum\limits_{\pi\in S_n} q^{\maj(\pi)} t^{des(\pi^{-1})}. \leqno(1)
$$
$$
\sum\limits_{\pi\in S_n} q^{\ell(\pi)} t^{\maj(\pi^{-1})}=
\sum\limits_{\pi\in S_n} q^{\maj(\pi)} t^{\maj(\pi^{-1})} . \leqno(2)
$$
\end{cor}


\bigskip



An alternative proof of
Theorem~\ref{FS} may be obtained using
the following classical fact~\cite[Prop. 1.3.17]{ECI}.

\begin{fac}\label{Fact3}
For any $M\subseteq [n-1]$,
$$
\sum_{\{\pi\in S_n\,|\,\Des(\pi^{-1})\subseteq M\}} q^{inv(\pi)}
= {n \brack \Delta M}_q. 
$$
\end{fac}

Garsia and Gessel proved that a similar identity holds for the major index.

\begin{thm}\label{GG}{\rm \cite[Theorem 3.1]{GaG}}
For any $M\subseteq [n-1]$,
$$
\sum_{\{\pi\in S_n\,|\,\Des(\pi^{-1})\subseteq M\}} q^{maj(\pi)}
= {n \brack \Delta M}_q. 
$$
\end{thm}

Combining this theorem with Fact~\ref{Fact3} implies Theorem~\ref{FS}.

\subsection{The Hyperoctahedral Group}

We denote by $B_{n}$ the group of all bijections $\sigma$ of the set 
$[\pm n]$ onto itself such that 
$$
\sigma(-a) = -\sigma(a) \qquad(\forall a\in[\pm n]),
$$
with composition as the group operation. This group is usually known as 
the group of ``signed permutations'' on $[n]$, or as the {\em hyperoctahedral 
group} of rank $n$. We identify $S_{n}$ as a subgroup of $B_{n}$, 
and $B_{n}$ as a subgroup of $S_{2n}$, in the natural ways.

If $\sigma \in B_{n}$ then write $\sigma = [a_{1}, \ldots ,a_{n}]$ 
to mean that $\sigma(i)=a_{i}$ for $1\le i\le n$, and let
\begin{eqnarray*}
\inv(\sigma) &:=& \inv (a_{1}, \ldots , a_{n} ),\\
\Des_A(\sigma)&:=& \Des(a_1,\ldots, a_n),\\
\des_{A}(\sigma) &:=& \des(a_{1}, \ldots , a_{n}),\\
\maj_A (\sigma) &:=& maj (a_{1}, \ldots , a_{n} ),\\
\Neg(\sigma) &:=& \{ i \in [n]: \; a_{i}<0 \},\\
neg(\sigma) &:=& |\Neg(\sigma)|.
\end{eqnarray*}
It is well known (see, e.g., \cite[Proposition 8.1.3]{BB}) that $B_{n}$
is a Coxeter group with respect to the generating set 
$\{ s_{0},s_{1},s_{2},\ldots,s_{n-1} \}$, where
\[ s_{0} := [-1,2, \ldots ,n ] \]
and
\[ s_{i} := [1,2,\ldots ,i-1,i+1, i,i+2, \ldots ,n] \qquad(1\le i < n). \]
This gives rise to two other natural statistics on $B_{n}$
(similarly definable for any Coxeter group), namely
\[ \ell_B(\sigma) := min \{ r \in \bbn: \; \sigma = s_{i_{1}} \ldots s_{i_{r}}
\;\; \mbox{ for some} \;\; i_{1}, \ldots ,i_{r} \in [0,n-1]
\} \]
(known as the {\em length} of $\sigma$) and
\[ \des_{B}(\sigma) := | \Des_B(\sigma)|, \]
where the {\it $B$-descent set} $\Des_B(\sigma)$ is defined as
\[ \Des_B(\sigma) := \{ i \in [0,n-1]\,|\,\ell_B(\sigma s_{i}) < \ell_B(\sigma) \}. \]

\begin{rem}\label{rem1}
Note that for every $\sg\in B_n$
$$
\Des_A(\sigma)= \Des_B(\sigma) \setminus \{0\}.
$$ 
\end{rem}

There are well known direct combinatorial ways to compute these statistics for 
$\sigma \in B_{n}$ (see, e.g., \cite[Propositions 8.1.1 and 8.1.2]{BB} 
or \cite[Proposition 3.1 and Corollary 3.2]{Bre94}), namely
\begin{equation}\label{e.length}
\ell_B(\sigma) = inv(\sigma) + \sum _{ i \in \Neg(\sigma)} |\sigma(i)|  
\end{equation}
and
\begin{equation}
\des_{B}(\sigma) = | \{ i \in [0,n-1]: \; \sigma(i) > \sigma(i+1) \} |,
\end{equation}
where  $\sigma(0) := 0$. 
For example, if $\sigma = [-3, 1,-6,2,-4,-5] \in B_{6}$ then 
$\inv(\sigma)=9$, $\des_{A}(\sigma)=3$, $\maj_A(\sigma)=11$, 
$neg(\sigma)=4$, $\ell_B(\sigma)=27$, and $\des_{B}(\sigma)=4$.

\section{Main Results}\label{main}

\begin{df}\label{d.fmaj}
The {\em flag major index} of a signed permutation $\sigma\in B_n$ is defined by
$$
fmaj(\sigma):=2\cdot \maj_A(\sigma)+neg(\sigma), 
$$
where $\maj_A(\sigma)$ is the major index of the sequence 
$(\sigma(1),\dots,\sigma(n))$ with respect to the order
$-n<\cdots<-1<1<\cdots<n$. 
\end{df}



\medskip

The main theorem is

\begin{thm}\label{t.sub_last}
For every subset $M\subseteq [0,n-1]$ and $i\in [\pm n]$
$$
\sum_{\{\sigma\in B_n\,|\,\Des_B(\sigma^{-1})\subseteq M,\,\sigma(n)=i\}}q^{\ell_B(\sigma)}=
\sum_{\{\sigma\in B_n\,|\,\Des_B(\sigma^{-1})\subseteq M,\,\sigma(n)=i\}}q^{\fmaj(\sigma)}.
$$
\end{thm}

\medskip

\noindent
See more details in Theorems~\ref{t.5.1} and~\ref{t.5.2} below. 
Forgetting $\sigma(n)$, we have

\begin{thm}\label{t.subFSB}
For every subset $M=\{m_1,\dots,m_t\}_< \subseteq [0,n-1]$,
\begin{eqnarray*}
\sum_{\{\sigma\in B_n\,|\,\Des_B(\sigma^{-1})\subseteq M\}}q^{\ell_B(\sigma)}
&=&
\sum_{\{\sigma\in B_n\,|\,\Des_B(\sigma^{-1})\subseteq M\}}q^{fmaj(\sigma)} =\\
&=&
{n \brack \Delta M}_q \cdot \prod\limits_{j=m_1+1}^n (1+q^j).
\end{eqnarray*}
\end{thm}

\medskip

We deduce a Foata-Sch\"utzenberger type theorem for $B_n$.

\medskip

\begin{thm}\label{t.FSB}
For every subset $M\subseteq [0,n-1]$,
$$
\sum_{\{\sigma\in B_n\,|\,\Des_B(\sigma^{-1})=M\}}q^{\ell_B(\sigma)}=
\sum_{\{\sigma\in B_n\,|\,\Des_B(\sigma^{-1})= M\}}q^{\fmaj(\sigma)} .
$$
\end{thm}

\section{Proofs}\label{s.proofs}


\begin{lem}\label{t.ObsB2}\
Let $\sigma\in B_n$ and $M=\{m_1,\dots,m_t\}_< \subseteq [0,n-1]$. 
Denote $m_{t+1} := n$.
Then: 
$\Des_B(\sigma^{-1})\subseteq M$ if and only if 
there exist (unique) integers $r_1,\ldots,r_t$ satisfying
$m_i\le r_i\le m_{i+1}$ $(\forall i)$ 
such that $(\sigma(1),\ldots,\sigma(n))$ is a shuffle of 
the following increasing sequences:
$$
(1,2,\dots,m_1),
$$ 
$$
(-r_1,-r_1+1,\dots,-(m_1+1)),
$$
$$
(r_1+1,r_1+2,\dots,m_2),
$$ 
$$
\vdots
$$
$$
(-r_t,-r_t+1,\dots,-(m_t+1)),
$$
$$
(r_t+1,r_t+2,\dots,n\,(\,=\,m_{t+1}\,)).
$$
Some of these sequences may be empty (if $r_i = m_i$ or $r_i = m_{i+1}$
for some $i$, or if $m_1=0$).
\end{lem}

\noindent
{\bf Proof.} 
Assume first that $0\in M$ (i.e., $m_1=0$). Then
$$
\Des_B(\sigma^{-1})\subseteq M \iff 
\Des_A(\sigma^{-1})\subseteq \{m_2,\ldots,m_t\}.
$$
This is equivalent to 
$$
\sigma^{-1}(j)<\sigma^{-1}(j+1)\qquad(\forall j\in [n]\setminus M),
$$
namely,
$$
\sigma^{-1}(m_i+1)<\sigma^{-1}(m_i+2)<\ldots<\sigma^{-1}(m_{i+1})
\qquad(1\le i\le t).
$$
Defining $m_i\le r_i\le m_{i+1}$ such that $\sigma^{-1}(r_i)$ is 
the last negative value in this sequence, 
we get the desired conclusion (with an empty sequence $(1,\ldots,m_1)$).

If $0\not\in M$ (i.e., $m_1>0$) then we get the same conclusion
for the intervals $[m_i+1,m_{i+1}]$ with $1\le i\le t$, and in addition
$\sigma^{-1}(1)>0$ so that
$$
0<\sigma^{-1}(1)<\sigma^{-1}(2)<\ldots<\sigma^{-1}(m_1).
$$
This gives us the desired conclusion (with a nonempty sequence $(1,\ldots,m_1)$).

\qed

\bigskip

\begin{thm}\label{t.5.1}
Let $M = \{ m_1,m_2,\ldots,m_t \}_{<} \subseteq [0,n-1]$ and $i \in [\pm n]$. 
Denote $m_0 := 0$ and $m_{t+1} := n$. Then
$$
\sum_{\{ \sigma\in B_{n}\,|\,Des_{B}(\sigma^{-1})\subseteq M, \;\sigma(n)=i \}}
q^{fmaj(\sigma)} =
\frac{\alpha_i(M)}{q^{n}-q^{-n}} \cdot 
{n \brack \Delta M}_q \cdot \prod_{j=m_1+1}^{n} (1+q^j), 
$$
where {\rm
$$
\alpha_i(M) = \cases{%
q^{m_1} - q^{-m_1},&if $i=m_1 >0$;\cr
q^{-m_s} - q^{-m_{s+1}},&if $i=m_{s+1}$ for $s\in[t]$;\cr
q^{m_{s+1}} - q^{m_s},&if $i=-(m_s+1)$ for $s\in[t]$;\cr
0, &otherwise.} 
$$}%
Note that the first case ($i=m_1 >0$) occurs only if $0\not\in M$.
\end{thm}

\bigskip

\noindent
{\bf Proof.} 
By induction on $n$. It is easy to verify the result for $n\le 2$.
Assume that $n\ge 3$, and that the result holds for $n-1$.

We shall use the induction hypothesis by ``deleting'' the value $\sigma(n)=i$
from $\sigma\in B_n$. Formally, for $i\in[\pm n]$, define a function
$\phi_i:[-n,n]\to[-(n-1),n-1]$ by
$$
a' = \phi_i(a) :=\cases{%
a, &if $|a|<|i|$;\cr
\frac{a}{|a|}(|a|-1), &if $|a|\ge |i|$ 
}\qquad(\forall a\in[-n,n]).
$$
When restricted to $[-n,n]\setminus\{0,i,-i\} = [\pm n]\setminus\{i,-i\}$,
$\phi_i$ is a bijection onto $[\pm (n-1)]$.
For $\sigma\in B_n$ with $\sigma(n)=i$, let
$$
\tau := [{\sigma(1)}',\ldots,{\sigma(n-1)}']\in B_{n-1}.
$$
Lemma~\ref{t.ObsB2} (for $\sigma$ and $\tau$) implies that, 
for any $M = \{m_1,\ldots,m_t\}_< \subseteq[0,n-1]$, 
the map $\sigma\mapsto\tau$ is a bijection 
from $\{\sigma\in B_n\,|\,\Des_B(\sigma^{-1})\subseteq M\}$
onto $\{\tau\in B_{n-1}\,|\,\Des_B(\tau^{-1})\subseteq M'\}$,
where $M' := \{m'_1,\ldots,m'_t\}_{\le}$.
Note that we may have $m'_s=m'_{s+1}$ (if $m_s+1=i=m_{s+1}$),
but this will make no difference in the sequel.

Assume that $\sigma(n-1)=j$, so that $\tau(n-1)=j'$. Then:
$$
\fmaj(\sigma) = \fmaj(\tau) + 2(n-1)\chi(j>i) + \chi(i<0).
$$
Denote
$$
B_n(M,i) := \{\sigma\in B_n\,|\,\Des_B(\sigma^{-1})\subseteq M,\,\sigma(n)=i\},
$$
and similarly $B_{n-1}(M',j')$.
Denote also
$$
M_{\pm} := \{m_1,\ldots,m_{t+1}\} \cup \{-(m_1+1),\ldots,-(m_t+1)\}.
$$
By Lemma~\ref{t.ObsB2}, if $\Des_B(\sigma^{-1})\subseteq M$ 
then $\sigma(n)\in M_{\pm}$. Thus the sum in the statement of our theorem
is zero whenever $i\not\in M_{\pm}$. We shall check the values $i\in M_{\pm}$
case-by-case.

First note that the coefficients $\alpha_i(M)$ defined in the statement of the 
theorem satisfy
\begin{equation}\label{e.alpha1}
\sum_{j=-n}^{n} \alpha_j(M) = q^n - q^{-n}
\end{equation}
and
\begin{equation}\label{e.alpha2}
\sum_{j>i} \alpha_j(M) = \cases{%
q^{-m_s} - q^{-n},& if $i=m_s$ for $s\in [t+1]$;\cr
q^{m_s} - q^{-n},& if $i=-(m_s+1)$ for $s\in [t]$.}
\end{equation}

\medskip

\noindent
{\bf Case 1:} $i=m_{s+1}$ for $s\in [t]$.

\noindent
Here
$$
\fmaj(\sigma) = \fmaj(\tau) + 2(n-1)\chi(j>i),
$$
so that
$$
\sum_{\sigma\in B_n(M,i)} q^{\fmaj(\sigma)} = 
\sum_{j\in [\pm n]\setminus\{i,-i\}} \,
\sum_{\sigma\in B_n(M,i)\atop \sigma(n-1)=j} q^{\fmaj(\sigma)} = 
$$
$$
= \sum_{j<i} \,\sum_{\tau\in B_{n-1}(M',j')} q^{\fmaj(\tau)}
+ \sum_{j>i} \,\sum_{\tau\in B_{n-1}(M',j')} q^{\fmaj(\tau)+2(n-1)} =
$$
$$
= \frac{{n-1 \brack \Delta M'}_q}{q^{n-1}-q^{-(n-1)}} \cdot  
\prod_{j=m'_1+1}^{n-1} (1+q^j) \cdot 
    \left[\sum_{j'\le i'} \alpha_{j'}(M') + 
q^{2(n-1)}\sum_{j'>i'} \alpha_{j'}(M')\right].
$$
We used here the fact that $j<i \iff j'\le i'$.
Now, by equalities (\ref{e.alpha1}) and (\ref{e.alpha2}) 
(for $M'$ instead of $M$, with $i'=m'_{s+1}=m_{s+1}-1$),
$$
\sum_{j'\le i'} \alpha_{j'}(M') + q^{2(n-1)}\sum_{j'>i'} \alpha_{j'}(M') =
$$
$$
= \left(q^{n-1}-q^{-m'_{s+1}}\right) 
+ q^{2(n-1)}\left(q^{-m'_{s+1}}-q^{-(n-1)}\right) =
$$
$$
= \left(q^{2(n-1)}-1\right)q^{-(m_{s+1}-1)}.
$$
Thus (using $m'_1=m_1$)
$$
\sum_{\sigma\in B_n(M,i)} q^{\fmaj(\sigma)} =
\frac{q^{2(n-1)}-1}{q^{n-1}-q^{-(n-1)}} \cdot q^{-(m_{s+1}-1)} 
\cdot {n-1 \brack \Delta M'}_q \cdot \prod_{j=m_1+1}^{n-1} (1+q^j) =
$$
$$
= q^{(n-1)-(m_{s+1}-1)} 
\cdot {n \brack \Delta M}_q \cdot \frac{[m_{s+1}-m_s]_q}{[n]_q}
\cdot \prod_{j=m_1+1}^{n} (1+q^j) \cdot \frac{1}{1+q^n} =
$$
$$
= q^{n-m_{s+1}} \cdot \frac{1-q^{m_{s+1}-m_s}}{1-q^{2n}}
\cdot {n \brack \Delta M}_q 
\cdot \prod_{j=m_1+1}^{n} (1+q^j) =
$$
$$
= \frac{q^{-m_s}-q^{-m_{s+1}}}{q^n-q^{-n}}
\cdot {n \brack \Delta M}_q \cdot \prod_{j=m_1+1}^{n} (1+q^j),
$$
as claimed.

\medskip

\noindent
{\bf Case 2:} $i=m_1 > 0$.

\noindent
The computations are as in the previous case, except that $m'_1=m_1-1$:
$$
\sum_{\sigma\in B_n(M,i)} q^{\fmaj(\sigma)} =
\frac{q^{2(n-1)}-1}{q^{n-1}-q^{-(n-1)}} \cdot q^{-(m_{1}-1)} 
\cdot {n-1 \brack \Delta M'}_q \cdot \prod_{j=m_1}^{n-1} (1+q^j) =
$$
$$
= q^{(n-1)-(m_{1}-1)} 
\cdot {n \brack \Delta M}_q \cdot \frac{[m_{1}]_q}{[n]_q}
\cdot \prod_{j=m_1+1}^{n} (1+q^j) \cdot \frac{1+q^{m_1}}{1+q^n} =
$$
$$
= q^{n-m_{1}} \cdot \frac{1-q^{2m_1}}{1-q^{2n}}
\cdot {n \brack \Delta M}_q 
\cdot \prod_{j=m_1+1}^{n} (1+q^j) =
$$
$$
= \frac{q^{m_1}-q^{-m_1}}{q^n-q^{-n}}
\cdot {n \brack \Delta M}_q \cdot \prod_{j=m_1+1}^{n} (1+q^j),
$$
as claimed.

\medskip

\noindent
{\bf Case 3:} $i=-(m_s+1)$ for $s\in [t]$.

\noindent
Here $\chi(i<0)=1$, so that
$$
\sum_{\sigma\in B_n(M,i)} q^{\fmaj(\sigma)} =
$$
$$
= \sum_{j<i} \,\sum_{\tau\in B_{n-1}(M',j')} q^{\fmaj(\tau)+1}
+ \sum_{j>i} \,\sum_{\tau\in B_{n-1}(M',j')} q^{\fmaj(\tau)+2(n-1)+1} =
$$
$$
= \frac{q \cdot {n-1 \brack \Delta M'}_q}{q^{n-1}-q^{-(n-1)}} \cdot  
\prod_{j=m'_1+1}^{n-1} (1+q^j) \cdot 
    \left[\sum_{j'\le i'} \alpha_{j'}(M') + 
q^{2(n-1)}\sum_{j'>i'} \alpha_{j'}(M')\right].
$$
By equalities (\ref{e.alpha1}) and (\ref{e.alpha2}) 
(for $M'$ instead of $M$, with $i'=-(m_s+1)'=-m_s$),
$$
\sum_{j'\le i'} \alpha_{j'}(M') + q^{2(n-1)}\sum_{j'>i'} \alpha_{j'}(M') =
$$
$$
= \left(q^{n-1}-q^{m'_s}\right) 
+ q^{2(n-1)}\left(q^{m'_s}-q^{-(n-1)}\right) =
$$
$$
= \left(q^{2(n-1)}-1\right)q^{m_s}.
$$
Thus (using $m'_1=m_1$)
$$
\sum_{\sigma\in B_n(M,i)} q^{\fmaj(\sigma)} =
\frac{q\left(q^{2(n-1)}-1\right)}{q^{n-1}-q^{-(n-1)}} \cdot q^{m_s} 
\cdot {n-1 \brack \Delta M'}_q \cdot \prod_{j=m_1+1}^{n-1} (1+q^j) =
$$
$$
= q^{n+m_s} 
\cdot {n \brack \Delta M}_q \cdot \frac{[m_{s+1}-m_s]_q}{[n]_q}
\cdot \prod_{j=m_1+1}^{n} (1+q^j) \cdot \frac{1}{1+q^n} =
$$
$$
= q^{n+m_s} \cdot \frac{1-q^{m_{s+1}-m_s}}{1-q^{2n}}
\cdot {n \brack \Delta M}_q 
\cdot \prod_{j=m_1+1}^{n} (1+q^j) =
$$
$$
= \frac{q^{m_{s+1}}-q^{m_s}}{q^n-q^{-n}}
\cdot {n \brack \Delta M}_q \cdot \prod_{j=m_1+1}^{n} (1+q^j),
$$
as claimed.

\qed

To prove the corresponding result for $\ell_{B}$ we will find it useful to use 
the fact that
\begin{equation}\label{e.ellB}
\ell_{B}(\sigma ) =\frac{\inv(\overline{\sigma})+neg(\sigma )}{2}
\end{equation}
for all $\sigma \in B_{n}$, where 
$\overline{\sigma } \eqdef (\sigma (-n),\ldots, \sigma (-1),\sigma (1),
\ldots , \sigma (n))$. This formula, first observed by Incitti in \cite{Inc}, is 
easily seen to be equivalent to (\ref{e.length}). So, for example, 
if $\sigma =[-3,1,-6,2,-4,-5]$ 
then $inv(\overline{\sigma })=50$ 
and $\ell_{B}(\sigma ) = \frac{50+4}{2}=27$. 

\begin{thm}\label{t.5.2} 
Let $M = \{ m_1,m_2,\ldots,m_t \}_{<} \subseteq [0,n-1]$ and $i \in [\pm n]$. 
Then $q^{\ell_B(\sigma)}$ satisfies exactly the same formula 
as does $q^{fmaj(\sigma)}$ in Theorem~\ref{t.5.1}.
\end{thm}

\noindent{\bf Proof.} 
By induction on $n$. As in the previous proof, the result is easy to verify 
for $n\le 2$. We assume that $n\ge 3$, and that the result holds for $n-1$.
Again, we may assume that 
$$
i\in M_{\pm} := \{m_1,\ldots,m_{t+1}\} \cup \{-(m_1+1),\ldots,-(m_t+1)\}.
$$
These values will be checked case-by-case.

\medskip

\noindent
{\bf Case 1:} $i=m_{s+1}$ for $s\in [t]$.

\noindent
By (\ref{e.ellB})
$$
\ell_B(\sigma) =\ell_B(\tau) + \frac{2(n-i)+0}{2}, 
$$ 
so that
$$
\sum_{\sigma\in B_n(M,i)} q^{\ell_B(\sigma)} = 
\sum_{j\in [\pm n]\setminus\{i,-i\}} \,
\sum_{\sigma\in B_n(M,i)\atop \sigma(n-1)=j} q^{\ell_B(\sigma)} = 
$$
$$
= \sum_{j\ne \pm i} \,\sum_{\tau\in B_{n-1}(M',j')} q^{\ell_B(\tau)+(n-i)} =
$$
$$
= \frac{q^{n-i} \cdot {n-1 \brack \Delta M'}_q}{q^{n-1}-q^{-(n-1)}} \cdot  
\prod_{j=m'_1+1}^{n-1} (1+q^j) \cdot 
\sum_{j'=-(n-1)}^{n-1} \alpha_{j'}(M') =
$$
$$
= q^{n-m_{s+1}} 
\cdot {n-1 \brack \Delta M'}_q \cdot \prod_{j=m_1+1}^{n-1} (1+q^j) =
$$
$$
= q^{n-m_{s+1}} 
\cdot {n \brack \Delta M}_q \cdot \frac{[m_{s+1}-m_s]_q}{[n]_q}
\cdot \prod_{j=m_1+1}^{n} (1+q^j) \cdot \frac{1}{1+q^n} =
$$
$$
= \frac{q^{-m_s}-q^{-m_{s+1}}}{q^n-q^{-n}}
\cdot {n \brack \Delta M}_q \cdot \prod_{j=m_1+1}^{n} (1+q^j),
$$
as claimed.

\medskip

\noindent
{\bf Case 2:} $i=m_1 > 0$.

\noindent
The computations are as in the previous case, except that $m'_1=m_1-1$:
$$
\sum_{\sigma\in B_n(M,i)} q^{\ell_B(\sigma)} = q^{n-m_1} 
\cdot {n-1 \brack \Delta M'}_q \cdot \prod_{j=m_1}^{n-1} (1+q^j) =
$$
$$
= q^{n-m_1} 
\cdot {n \brack \Delta M}_q \cdot \frac{[m_1]_q}{[n]_q}
\cdot \prod_{j=m_1+1}^{n} (1+q^j) \cdot \frac{1+q^{m_1}}{1+q^n} =
$$
$$
= \frac{q^{m_1}-q^{-m_1}}{q^n-q^{-n}}
\cdot {n \brack \Delta M}_q \cdot \prod_{j=m_1+1}^{n} (1+q^j),
$$
as claimed.

\medskip

\noindent
{\bf Case 3:} $i=-(m_s+1)$ for $s\in [t]$.

\noindent
Here 
$$
\ell_B(\sigma) = \ell_B(\tau) + \frac{(4(|i|-1)+2(n-|i|)+1)+1}{2}.
$$ 
Thus, using $m'_1=m_1$:
$$
\sum_{\sigma\in B_n(M,i)} q^{\ell_B(\sigma)} =
\sum_{j\ne \pm i} \,\sum_{\tau\in B_{n-1}(M',j')} q^{\ell_B(\tau)+(n+|i|-1)} =
$$
$$
= \frac{q^{n+|i|-1} \cdot {n-1 \brack \Delta M'}_q}{q^{n-1}-q^{-(n-1)}} \cdot  
\prod_{j=m'_1+1}^{n-1} (1+q^j) \cdot \sum_{j'=-(n-1)}^{n-1} \alpha_{j'}(M') =
$$
$$
= q^{n+m_s} 
\cdot {n-1 \brack \Delta M'}_q \cdot \prod_{j=m_1+1}^{n-1} (1+q^j) =
$$
$$
= q^{n+m_s} 
\cdot {n \brack \Delta M}_q \cdot \frac{[m_{s+1}-m_s]_q}{[n]_q}
\cdot \prod_{j=m_1+1}^{n} (1+q^j) \cdot \frac{1}{1+q^n} =
$$
$$
= \frac{q^{m_{s+1}}-q^{m_s}}{q^n-q^{-n}}
\cdot {n \brack \Delta M}_q \cdot \prod_{j=m_1+1}^{n} (1+q^j),
$$
as claimed.

\qed

\medskip

From Theorems \ref{t.5.1} and \ref{t.5.2} 
we immediately deduce Therem~\ref{t.sub_last}.



By summing Theorems \ref{t.5.1} and \ref{t.5.2} over $i \in [\pm n]$
we obtain Theorem~\ref{t.subFSB}, 
which was the original motivation for this work, and which is the analogue, 
for the hyperoctahedral group, of Theorem~\ref{GG}.


It would be interesting to have combinatorial (bijective) proofs of these results. 

\medskip

\noindent
{\bf Added in Proof:}
A bijective proof of Theorem~\ref{t.FSB}
has been found by Foata and Han~\cite{FH3}.

\section{Final Remarks}

\subsection{Open Problems}

Numerical evidence suggests that the following holds.

\begin{con}
For every subset $M \subseteq [0,n-1]$ and $i\in [ \pm n]$ the polynomial 
$$
\sum_{\{\sigma\in B_n\,|\,\Des_B(\sigma^{-1}) \subseteq M,\,\sigma(n)=i\}}q^{\ell_B(\sigma)}=
\sum_{\{\sigma\in B_n\,|\,\Des_B(\sigma^{-1}) \subseteq M,\,\sigma(n)=i\}}q^{\fmaj(\sigma)}.
$$
is (symmetric and) unimodal.
\end{con}

\begin{con}
For every subset $M \subseteq [0,n-1]$ the polynomial 
$$
\sum_{\{\sigma\in B_n\,|\,\Des_B(\sigma^{-1}) \subseteq M \}}q^{\ell_B(\sigma)}=
\sum_{\{\sigma\in B_n\,|\,\Des_B(\sigma^{-1}) \subseteq M \}}q^{\fmaj(\sigma)}.
$$
is (symmetric and) unimodal.
\end{con}

Using well known results (see, e.g., \cite[Proposition 1 and Theorem 11]{Sta}) 
and Theorems \ref{t.subFSB}, \ref{t.5.1} and \ref{t.5.2}, 
it is easy to see that the above conjectures are equivalent to the following.

\begin{con}\label{c.f2} 
For $0 \leq k < n$ the polynomial 
$$
{n\brack k}_q  \cdot \prod\limits_{j=k+1}^n (1+q^j) 
$$
is (symmetric and) unimodal. 
\end{con}

We have verified these conjectures for $n \leq 15$. Conjecture \ref{c.f2} 
clearly holds for $k=n-1$ and is known to be true for $k=0$ (see, e.g., 
\cite[p. 510]{Sta}). 
We have checked that the polynomials 
$$
\sum_{\{\sigma\in B_n\,|\,\Des_B(\sigma^{-1}) = M,\,\sigma(n)=i\}}q^{\ell_B(\sigma)}=
\sum_{\{\sigma\in B_n\,|\,\Des_B(\sigma^{-1}) = M,\,\sigma(n)=i\}}q^{\fmaj(\sigma)}.
$$
and
$$
\sum_{\{\sigma\in B_n\,|\,\Des_B(\sigma^{-1}) = M \}}q^{\ell_B(\sigma)}=
\sum_{\{\sigma\in B_n\,|\,\Des_B(\sigma^{-1}) = M \}}q^{\fmaj(\sigma)}.
$$
are unimodal for all $M \subseteq [0,n-1]$ and $i\in [ \pm n]$ if 
$n \leq 5$. In general, these polynomials are not symmetric.

\subsection{Classical Weyl Groups of Type $D$}


Let $D_n$ be the classical Weyl group of type $D$ and rank $n$.
For an element $\sg\in D_n$, let $\ell_D(\sg)$ 
be the length of $\sg$ with respect to the Coxeter generators of $D_n$.
It is well known that we may take
$$
D_n=\{\sg\in B_n\,|\,neg(\sg)\equiv 0 \hbox{ mod }2\}.
$$
Let $\sg=[\sg(1),\dots,\sg(n)]\in D_n$. Biagioli and Caselli~\cite{BC} 
introduced a flag major index for $D_n$
$$
\fmaj_D(\sg):=\fmaj(\sg(1),\dots,\sg(n-1),|\sg(n)|).
$$
By definition,
\begin{equation}\label{e.fmajD}
\sum_{\sigma\in D_n}q^{fmaj_D(\sigma)} = 
\sum_{\{\sigma\in B_n\,|\,\sigma(n)>0\}}q^{fmaj(\sigma)}.
\end{equation}

\begin{pro}\label{t.d1}
$$
\sum_{\sigma\in D_n}q^{\ell_D(\sigma)} =
\sum_{\{\sigma\in B_n\,|\,\sigma(n)>0\}}q^{\ell_B(\sigma)}.
$$
\end{pro}


\noindent{\bf Proof.}
It is well known (see, e.g.~\cite[\S 3.15]{Hu}) that 
$$
\sum\limits_{\sigma\in D_n}q^{\ell_D(\sigma)}=
[n]_q \cdot \prod\limits_{i=1}^{n-1}[2i]_q.
$$
On the other hand, 
$$
\sum_{\{\sigma\in B_n\,|\,\sigma(n)>0\}}q^{\ell_B(\sigma)}=
\sum_{i=1}^n\,\sum_{\{\sigma\in B_n\,|\,\sigma(n)= i\}}q^{\ell_B(\sigma)} =
$$
$$
= \sum_{i=1}^n\, 
\sum_{\{\sigma\in B_n\,|\,\Des_B(\sigma^{-1})\subseteq [0,n-1],\,\sigma(n)= i\}}q^{\ell_B(\sigma)}.
$$
By Theorem~\ref{t.5.2}, this is equal to
$$
\frac{q^0-q^{-n}}{q^n-q^{-n}} \cdot
{n \brack 1,\ldots,1}_q \cdot \prod_{j=1}^{n} (1+q^j)
= [n]_q ! \cdot \prod_{j=1}^{n-1} (1+q^j)
= [n]_q \cdot \prod_{j=1}^{n-1} [2j]_q,
$$
completing the proof.

\qed

\medskip

We deduce the following type $D$ analogue (first proved in \cite{BC})
of MacMahon's theorem.

\begin{cor}
$$
\sum\limits_{\sigma\in D_n}q^{\fmaj_D(\sigma)}
= \sum\limits_{\sigma\in D_n}q^{\ell_D(\sigma)}.
$$
\end{cor}

\noindent{\bf Proof.}
Combine (\ref{e.fmajD}) and Proposition~\ref{t.d1} with Theorem~\ref{t.sub_last}.

\qed

\begin{prb}
Find an analogue of the Foata-Sch\"utzenberger theorem for $D_n$.
\end{prb}

The obvious candidate for such an analogue does not work.




\subsection{Two Versions of the Flag Major Index}

The flag major index of $\sg\in B_n$, $\hbox{\it flag-major}(\sigma)$, 
was originally defined as 
the length of a distinguished canonical expression for $\sigma$. 
In~\cite{AR} this length was shown to be equal to
$2\cdot \maj_A(\sigma)+neg(\sigma)$,
where the major index of the sequence $(\sigma(1),\ldots,\sigma(n))$ 
was taken with respect to the order
$-1<\cdots<-n<1<\cdots<n$.
In~\cite{ABR} we considered a different order:
$-n<\cdots<-1<1<\cdots<n$
(i.e., we defined $fmaj$ as in Section~\ref{main} above). 

While both versions give type $B$ analogues of the 
MacMahon and Carlitz identities, only the second one gives
an analogue of the Foata-Sch\"utzenberger theorem.
On the other hand, the first one has the alternative natural 
interpretation as length, as mentioned above,
and also produces a natural analogue of the signed
Mahonian formula of Gessel and Simion, see~\cite{AGR}.
The relation between these two versions and their (possibly different)
algebraic roles requires further study.


\section{Appendix}

In this appendix we give an alternative proof of Theorems~\ref{t.subFSB} 
and~\ref{t.FSB} (but not Theorem~\ref{t.sub_last}), 
using $q$-binomial identities.

\subsection{Binomial Identities}

In this subsection we recall several $q$-binomial identities.

\medskip

\begin{lem}\label{t.q2}
For every positive integer $n$
$$
\sum_{k=0}^n {n\brack k}_{q^2} q^k = \prod_{i=1}^n (1+q^i). 
$$
\end{lem}

This identity may be easily proved by induction on $n$.
The following lemma is a multinomial extension of it.

\begin{lem}\label{t.mq2}
For every subset $M=\{m_1,\dots,m_t\}_<\subseteq [0,n-1]$
$$
\sum\limits_{r_1,\ldots,r_t} {n\brack \Delta M_r}_{q^2}
q^{\sum_{i=1}^{t} (r_i-m_i)} =
{n \brack \Delta M}_q \cdot \prod_{j=m_1+1}^n (1+q^j),
$$
where $m_{t+1}:=n$, the sum on the left hand side is over 
all $r_1,\ldots,r_t$ such that $m_i\le r_i\le m_{i+1}$ $(\forall i)$,
and
$$
{n\brack \Delta M_r}_{q^2} := 
{n\brack m_1,r_1-m_1,m_2-r_1,\ldots,r_t-m_t,m_{t+1}-r_t}_{q^2}.
$$
\end{lem}

\noindent{\bf Proof.} 
%
Decomposing the multinomial coefficient,
$$
\sum\limits_{r_1,\ldots,r_t} 
{n\brack \Delta M_r}_{q^2} q^{\sum_{i=1}^{t} (r_i-m_i)}
= {n\brack \Delta M}_{q^2} \cdot \prod_{i=1}^{t} \sum_{r_i=m_i}^{m_{i+1}} 
{m_{t+1}-m_t \brack r_t-m_t}_{q^2} q^{r_i-m_i}.
$$
By Lemma~\ref{t.q2} this is equal to
$$
{n \brack \Delta M}_{q^2} \cdot \prod_{i=1}^{t} \prod_{j=1}^{m_{i+1}-m_i} (1+q^j)
$$
and, since
$$
{n \brack \Delta M}_{q^2} = {n \brack \Delta M}_q \cdot 
\prod_{j=1}^{n} (1+q^j) \cdot 
\left[ \prod_{i=0}^{t} \prod_{j=1}^{m_{i+1}-m_i} (1+q^j)\right]^{-1} =
$$
$$
= {n \brack \Delta M}_q \cdot 
\prod_{j=m_1+1}^{n} (1+q^j) \cdot 
\left[ \prod_{i=1}^{t} \prod_{j=1}^{m_{i+1}-m_i} (1+q^j)\right]^{-1},
$$
we get the desired conclusion. 

\qed

\bigskip

The following ``$q$-binomial Theorem'' is well known.

\begin{thm}\label{t.q-bin}
$$
\prod\limits_{i=1}^n (1+q^i x)= \sum\limits_{k=0}^n {n\brack k}_q q^{k+1\choose 2} x^k.
$$
\end{thm}

\subsection{An Alternative Proof of Theorems~\ref{t.subFSB} and~\ref{t.FSB}} 


\noindent{\bf Proof of Theorem~\ref{t.subFSB}.}
Let $m_0:=0$ and $m_{t+1}:=n$.
By Lemma~\ref{t.ObsB2}, for each $\sigma\in B_n$ with
$\Des_B(\sigma^{-1})\subseteq M$ 
there exist $r_1,\ldots,r_t$ such that $m_i\le r_i\le m_{i+1}$ $(\forall i)$
and $(\sigma(1),\ldots,\sigma(n))$ is a shuffle of 
the following increasing sequences:
$$
(1,\dots,m_1),
$$
$$
(-r_1,\dots,-(m_1+1)),
$$
$$
(r_1+1,\dots,m_2),
$$
$$
\vdots
$$
$$
(-r_t,\dots,-(m_t+1)),
$$
$$
(r_t+1,\dots,n).
$$
By~(\ref{e.length}),
$$
\sum_{\{\sigma\in B_n\,|\,\Des_B(\sigma^{-1})\subseteq M\}}q^{\ell_B(\sigma)}=
\sum_{\{\sigma\in B_n\,|\,\Des_B(\sigma^{-1})\subseteq M\}}q^{\inv(\sigma)+
\sum_{\sg(i)<0} |\sg(i)|}.
$$
Note that
$$
\sigma(i)<0 \iff (\exists j)\; m_j+1\le |\sigma(i)|\le r_j.
$$
Therefore
$$
{\sum_{\sg(i)<0} |\sg(i)|} =
{\sum_{i=1}^t [(m_i+1) + \ldots + r_i]} =
{\sum_{i=1}^t {1\over 2}(r_i-m_i)(r_i+m_i+1)}.
$$
This is a constant, once we fix $r_1,\ldots,r_t$ (and $M$).
The inversion number of a shuffle does not depend on the actual values
of the elements in the shuffled sequences, but only on their order. 
Therefore, by Observation~\ref{ObsA} and Fact~\ref{Fact3},
$$
\sum\limits_{\sg} q^{\inv(\sg)}=
{n\brack m_1,r_1-m_1,m_2-r_1,\dots,r_t-m_t,m_{t+1}-r_t}_q,
$$
where the sum on the left hand side is over all $\sg\in B_n$
with $\Des_B(\sigma^{-1})\subseteq M=\{m_1,\dots m_t\}_<$
and prescribed $r_1,\ldots,r_t$.

Combining these two formulas, we get
$$
\sum\limits_{\{\sigma\in B_n\,|\,\Des_B(\sigma^{-1})\subseteq M\}}
q^{\ell_B(\sigma)}=
$$
$$
= \sum\limits_{r_1,\ldots,r_t} 
{n\brack m_1,r_1-m_1,m_2-r_1,\dots,m_{t+1}-r_t}_q 
q^{\sum_{i=1}^t {1\over 2}(r_i-m_i)(r_i+m_i+1)} =
$$
$$
= {n \brack \Delta M}_q \cdot
\prod\limits_{i=1}^t \sum_{r_i=m_i}^{m_{i+1}} 
{m_{i+1}-m_i\brack r_i-m_i}_q q^{{1\over 2}(r_i-m_i)(r_i+m_i+1)}.
$$
By the $q$-binomial Theorem (Theorem~\ref{t.q-bin}), 
with $x=q^{m_i}$ and $n=m_{i+1}-m_i$,
$$
\prod\limits_{j=m_i+1}^{m_{i+1}} (1+q^j)=
\sum\limits_{r_i=m_i}^{m_{i+1}} 
{m_{i+1}-m_i\brack r_i-m_i}_q q^{{r_i-m_i+1\choose 2}+m_i(r_i-m_i)} =
$$
$$
= \sum\limits_{r_i=m_i}^{m_{i+1}}
{m_{i+1}-m_i\brack r_i-m_i}_q q^{{1\over 2}(r_i-m_i)(r_i+m_i+1)}.
$$
Thus
$$
\sum_{\{\sigma\in B_n\,|\,\Des_B(\sigma^{-1})\subseteq M\}}q^{\ell_B(\sigma)} =
$$
$$
= {n \brack \Delta M}_q \cdot
\prod_{i=1}^{t} \prod\limits_{j=m_i+1}^{m_{i+1}} (1+q^j)
= {n \brack \Delta M}_q \cdot
\prod_{j=m_1+1}^{n} (1+q^j).
$$
This completes the proof of the second equality in the theorem, 
computing a generating function for $\ell_B$.

\bigskip

\noindent
An analogous computation holds for $\fmaj$: 
by Definition~\ref{d.fmaj}, Lemma~\ref{t.ObsB2} and Theorem~\ref{GG},
$$
\sum\limits_{\{\sigma\in B_n\,|\,\Des_B(\sigma^{-1})\subseteq M\}}
q^{\fmaj(\sigma)}=
\sum\limits_{\{\sigma\in B_n\,|\,\Des_B(\sigma^{-1})\subseteq M\}}
q^{2\cdot \maj(\sigma)+neg(\sigma)}=
$$
$$
= \sum\limits_{r_1,\ldots,r_t} 
{n\brack m_1,r_1-m_1,m_2-r_1,\ldots,r_t-m_t,m_{t+1}-r_t}_{q^2}
q^{\sum_{i=1}^t (r_i-m_i)}.
$$
By Lemma~\ref{t.mq2} this is equal to
$$
{n \brack \Delta M}_q \cdot \prod\limits_{i=m_1+1}^{n} (1+q^i),
$$
as claimed.

\qed

\bigskip

\noindent{\bf Proof of Theorem~\ref{t.FSB}.}
Apply the Inclusion-Exclusion Principle to Theorem~\ref{t.subFSB}. 

\qed

\end{document}